\documentclass[11pt]{amsart}
\usepackage{amsfonts}
\usepackage{amssymb,latexsym}
\usepackage{enumerate}
\usepackage{mathrsfs}
\usepackage{amsmath}
\usepackage{amsthm}

scaled\magstephalf

\newtheorem{Thm}{Theorem}[section]

\newtheorem{Lem}[Thm]{Lemma}
\newtheorem{Rem}[Thm]{Remark}

\numberwithin{equation}{section}

\begin{document}
\setlength{\baselineskip}{1.2\baselineskip}
\title[Gradient Estimates]{Gradient Estimates of Mean Curvature Equations with Neumann Boundary Condition}

\author{Xi-Nan Ma}
\address{Department of Mathematics\\
         University of Science and Technology of China\\
         Hefei Anhui 230026 CHINA}
         \email{xinan@ustc.edu.cn}
\author{Jinju Xu}
\address{Department of Mathematics\\
         University of Science and Technology of China\\
         Hefei Anhui 230026 CHINA}
         \email{july25@mail.ustc.edu.cn}

\thanks{Research of the first author was supported by NSFC  No.11125105 and Wu Wen-Tsun Key Laboratory of Mathematics. }

\maketitle

\begin{abstract}
In this paper, we use the maximum principle to get the gradient estimate for the solutions of the prescribed  mean curvature equation with Neumann boundary value problem, which gives a positive answer for the question raised by Lieberman \cite{Lieb13} in page 360.
As a consequence, we obtain the corresponding existence theorem for a class of mean curvature equations.
\end{abstract}

\section{Introduction}
Gradient estimate for the prescribed mean curvature equation has been extensively studied. The interior gradient estimate, for the minimal surface equation was obtained in the case of two variables by Finn \cite{FN54}. Bombieri, De Giorgi and M.Miranda \cite{BGM69}  obtained the estimate for high dimension case. For the general mean curvature equation, such estimate had also been obtained by Ladyzhenskaya and Ural'tseva \cite{LU70}, Trudinger \cite{TR73} and Simon \cite{Sim76}. All their methods were test function argument and a resulting Sobolev inequality. In 1983, Korevaar \cite{Kor86} introduced the normal variation technique and got the maximum principle proof for the interior gradient estimate on the minimal surface equation. Wang \cite{Wang98} gave a new proof for the mean curvature equation via  standard Bernstein technique.  The Dirichlet problem for the prescribed mean curvature equation had  been  studied by Jenkins-Serrin \cite{JS68} and Serrin \cite{S69}. More detailed history could be found in Gilbarg and Trudinger \cite{GT01}.

  For the  mean curvature equation with prescribed contact angle boundary value problem,  Ural'tseva \cite{Ur73} first got the boundary gradient estimates and the corresponding existence theorem. At the same time, Simon-Spruck \cite{SS76} and Gerhardt \cite{Ger76} also obtained existence theorem on the positive gravity case. For more general quasilinear divergence structure equation with conormal derivative boundary value problem, Lieberman \cite{Lie83}  gave the gradient estimate. They obtained these estimates also via test function technique.

  Spruck \cite{Sp75} used the maximum principle to obtain boundary gradient estimate in two dimension for the positive gravity capillary problems.  Korevaar \cite{Kor88} generalized his normal variation technique and got the gradient estimates for the positive gravity case in high dimension case. In \cite{Lieb84, Lieb87}, Lieberman developed the maximum principle approach on the boundary gradient estimates to the quasilinear elliptic equation with oblique derivative boundary value problem, and in \cite{Lieb88} he got the maximum principle proof for the gradient estimates on the  general quasilinear elliptic equation with capillary boundary value problems.

  In a recent book written by Lieberman (\cite{Lieb13}, in page 360), he posed the following question, how to get the gradient estimates for the mean curvature equation with Neumann boundary value problem. In this paper we use the technique developed by Spruck \cite{Sp75}, Lieberman \cite{Lieb88} and Wang \cite{Wang98} to get a positive answer. As a consequence, we obtain an existence theorem for a class of mean curvature equations with Neumann boundary value problem.

We first consider the boundary gradient estimates for the mean curvature equation with Neumann boundary value problem.
Now let's state our main gradient estimates.
\begin{Thm}\label{Thm1.1}
Suppose $u\in C^{2}(\overline\Omega)\bigcap C^{3}(\Omega)$ is a bounded  solution for the following boundary value problem
\begin{align}
 \texttt{div}(\frac{Du}{\sqrt{1+|Du|^2}}) =&f(x, u)   \quad\text{in}\quad \Omega, \label{1.1}\\
              \frac{\partial u}{\partial \gamma} = &\psi(x, u)  \quad\text{on} \quad\partial \Omega,\label{1.2}
\end{align}
where $\Omega \subset \mathbb R^n $ is a bounded domain, $n\geq 2$, $\partial \Omega \in C^{3}$, $\gamma$ is the inward unit normal to $\partial\Omega $.

We assume $f(x,z) \in C^{1}(\overline\Omega\times [-M_0, M_0])$ and $\psi(x,z) \in C^{3}(\overline\Omega\times [-M_0, M_0])$,
and there exist positive constants $M_0, L_1, L_2$ such that
\begin{align}
|u|\leq& M_0\quad in\quad\overline\Omega,\label{1.3}\\
f_z(x,z)\geq &0 \,\,\,\,\quad in\quad \overline\Omega\times[-M_0, M_0],\label{1.4}\\
|f(x,z)|+|f_{x}(x,z)|\leq & L_1 \quad in\quad \overline\Omega\times[-M_0, M_0],\label{1.5}\\
|\psi(x,z)|_{C^3(\overline\Omega\times[-M_0, M_0])}\leq&  L_2.\label{1.6}
\end{align}

Then there exists a small positive constant
$\mu_0$ such that we have the following estimate
$$\sup_{\overline\Omega_{\mu_0}}|Du|\leq \max\{M_1, M_2\},$$
where $M_1$ is a positive constant depending only on $n, \mu_0, M_0, L_1$, which is from the interior gradient estimates;
$M_2$ is  a positive constant depending only on $n, \Omega, \mu_0, M_0, L_1, L_2$, and $d(x) =\texttt{dist}(x, \partial\Omega), \Omega_{\mu_0} = \{x \in \Omega : d(x)<\mu_0\}.$
\end{Thm}
 As we stated before, there is a standard interior gradient estimates for the mean curvature equation.

\begin{Rem}[\cite{GT01}]\label{Rem1.1}
If $u\in C^{3}(\Omega)$ is a bounded solution for the equation \eqref{1.1}  with \eqref{1.3}, and if $f \in C^{1}(\overline\Omega \times [-M_0, M_0])$ satisfies the conditions \eqref{1.4}-\eqref{1.5}, then for any subdomain
$\Omega'\subset\subset\Omega$, we have
$$\sup_{\Omega'}|Du|\leq M_1,$$
where $M_1$ is a positive constant depending only on $n,  M_0, \texttt{dist} (\Omega', \partial\Omega), L_1$.
\end{Rem}

From the standard bounded estimates for the prescribed mean curvature equation in Concus-Finn \cite{CF74} ( see also Spruck \cite{Sp75}), we can get the following existence theorem for the Neumann boundary value problem of mean curvature equation.

\begin{Thm}\label{Thm1.2}
Let $\Omega \subset\mathbb R^n $ be a bounded domain, $n\geq 2$, $\partial \Omega\in C^{3}$, $\gamma$ is the inward unit normal to $\partial\Omega $.
If  $\psi \in C^{3}(\overline\Omega)$,  is a given function, then the following boundary value problem
\begin{align}
 \texttt{div}(\frac{Du}{\sqrt{1+|Du|^2}}) =& u   \quad\text{in}\quad \Omega, \label{1.7}\\
              \frac{\partial u}{\partial \gamma} = &\psi(x)  \quad\text{on} \quad\partial \Omega,\label{1.8}
\end{align}
 exists a unique solution $ u \in C^2(\overline\Omega)$.
\end{Thm}

The rest of the paper is organized as follows. In section 2, we first give the definitions and some notations. We prove the main Theorem~\ref{Thm1.1} in section 3 under the help of one lemma. In section 4, we give the proof of  Theorem~\ref{Thm1.2}.

\section{PRELIMINARIES}
We denote by $\Omega$ a bounded  domain in $\mathbb{R}^n$, $n\geq 2$,  $\partial \Omega\in C^{3}$,   set
\begin{align*}
 d(x)=\texttt{dist}(x,\partial \Omega),
 \end{align*}
 and
\begin{align*}
 \Omega_\mu=&\{{x\in\Omega:d(x)<\mu}\}.
 \end{align*}
Then it is well known that there exists a positive constant $\mu_{1}>0$ such that $d(x) \in C^3(\overline \Omega_{\mu_{1}})$. As in Simon-Spruck \cite{SS76} or Lieberman \cite{Lieb13} in page 331,  we can take $\gamma= D d$ in $\Omega_{\mu_{1}}$ and note that  $\gamma$ is a $C^2(\overline \Omega_{\mu_{1}})$ vector field. As mentioned in \cite{Lieb88} and the book \cite{Lieb13}, we also have the following formulas

\begin{align}\label{2.1}
\begin{split}
|D\gamma|+|D^2\gamma|\leq& C(n,\Omega) \quad\text{in}\quad \Omega_{\mu_{1}},\\
 \sum_{1\leq i\leq n}\gamma^iD_j\gamma^i=0,  \sum_{1\leq i\leq n}\gamma^iD_i\gamma^j=&0, \,|\gamma|=1 \quad\text{in} \quad\Omega_{\mu_{1}}.
\end{split}
\end{align}
As in \cite{Lieb13}, we define
 \begin{align}\label{2.2}
\begin{split}
c^{ij}=&\delta_{ij}-\gamma^i\gamma^j  \quad \text{in} \quad \Omega_{\mu_{1}},
\end{split}
\end{align}
 and for a vector $\zeta \in R^n$, we write $\zeta'$ for the vector with $i-$th component $ \sum_{1\leq j\leq n}c^{ij}\zeta_j$. So
 \begin{align}\label{2.3}
\begin{split}
|D'u|^2=& \sum_{1\leq i,j\leq n}c^{ij}u_iu_j.
\end{split}
\end{align}
Let
\begin{align}\label{2.4}
\begin{split}
a^{ij}(D u)=v^2\delta_{ij}-u_iu_j, \quad
v=(1+|D u|^2)^{\frac{1}{2}}.
\end{split}
\end{align}

Then the equations \eqref{1.1}, \eqref{1.2} are equivalent to the following boundary value problem
\begin{align}\label{2.5}
\sum_{i,j=1}^n a^{ij}u_{ij}=&f(x,u)v^3 \quad \text{in}\,\Omega,\\
u_{\gamma}=& \psi(x, u) \quad \text{on}\,\partial\Omega.\label{2.6}
\end{align}

\section{Proof of Theorem~\ref{Thm1.1} }

Now we begin to prove Theorem~\ref{Thm1.1}, as  mentioned in introduction, using the technique developed by Spruck \cite{Sp75}, Lieberman \cite{Lieb88} and Wang \cite{Wang98}. We shall choose an auxiliary function which contains  $|Du|^2$ and other lower order terms. Then we use the maximum principle for this auxiliary function in $\overline\Omega_{\mu_0}, 0<\mu_0<\mu_1$. At last,  we get our estimates.

{\em Proof of  Theorem~\ref{Thm1.1}.}

Setting $w=u-\psi(x,u)d$,  we choose the following auxiliary function
$$\Phi(x)=\log|Dw|^2e^{1+M_0+u}e^{\alpha_0 d}, \quad  x \in \overline \Omega_{\mu_{0}},$$ where $\alpha_0=|\psi|_{C^0(\overline\Omega\times[-M_0, M_0])}+C_0+1 $, $C_0$ is  a positive constant depending only on $n,\Omega$.

In order to simplify the computation, we let
\begin{align}\label{3.0}
 \varphi(x)= \log \Phi(x) =\log\log|Dw|^2+h(u)+g(d).
 \end{align}

In our case, we take
\begin{align}\label{3.1}
 h(u)=1+M_0+u,\quad g(d)=\alpha_0 d.
 \end{align}

 We assume that
$\varphi(x)$ attains its maximum at $x_0 \in \overline \Omega_{\mu_{0}}$, where $0<\mu_0<\mu_1$ is a sufficiently small number which we shall decide it  later.

Now we divide three cases to complete the proof of  Theorem~\ref{Thm1.1}.

Case I: If $\varphi(x)$ attains its maximum at $x_0 \in \partial\Omega$, then we shall use the Hopf Lemma to get the bound of $|Du|(x_0)$.

Case II: If $\varphi(x)$ attains its maximum at $x_0 \in\partial\Omega_{\mu_0}\bigcap\Omega$, then we shall get the estimates via the standard interior gradient bound \cite{GT01}.

Case III:  If $\varphi(x)$ attains its maximum at $x_0 \in \Omega_{\mu_0}$, in this case for the sufficiently small constant $\mu_0>0$,  then we can use the maximum principle to get the bound of $|Du|(x_0)$.

Now  all computations work at the point $x_0$.

{\bf Case 1.} If  $\varphi(x)$ attains its maximum at $x_0\in \partial \Omega$, we shall get the bound of $|Du|(x_0)$.

We differentiate $\varphi$ along the normal direction.

\begin{align}\label{3.2}
\begin{split}
\frac{\partial\varphi}{\partial\gamma}=&\frac{\sum_{1\leq i\leq n}(|Dw|^2)_i\gamma^i}{|Dw|^2\log|Dw|^2}+h'u_{\gamma}+g'.
\end{split}
\end{align}

Since
\begin{align}
w_i=&u_i-\psi_u u_i d -\psi_{x_i} d-\psi\gamma^i,\label{3.3}\\
|Dw|^2=&|D'w|^2+w^2_\gamma,\label{3.4}
\end{align}

we have
\begin{align}
w_\gamma=&u_\gamma-\psi_uu_\gamma d-\sum_{1\leq i\leq n}\psi_{x_i}\gamma^i d-\psi=0\quad\text{on}\quad\partial\Omega,\label{3.5}\\
(|Dw|^2)_i=&(|D'w|^2)_i\quad\text{on}\quad\partial\Omega.\label{3.6}
\end{align}

Applying  \eqref{2.1}, \eqref{2.3} and \eqref{3.6}, it follows that
\begin{align}\label{3.7}
\begin{split}
\sum_{1\leq i\leq n}(|Dw|^2)_i\gamma^i=&\sum_{1\leq i\leq n}(|D'w|^2)_i\gamma^i\\
=&2\sum_{1\leq i,k,l\leq n}c^{kl}w_{ki}w_l\gamma^i\\
=&2\sum_{1\leq i,k,l\leq n}c^{kl}u_{ki}u_l\gamma^i-2\sum_{1\leq k,l\leq n}c^{kl}u_lD_k\psi,
\end{split}
\end{align}
where
\begin{align*}
D_k\psi= \psi_{x_k} + \psi_u u_k.
\end{align*}

Differentiating \eqref{2.6} with respect to tangential direction,   we have
\begin{align}\label{3.8}
\sum_{1\leq k\leq n}c^{kl}(u_{\gamma})_k=&\sum_{1\leq k\leq n}c^{kl}D_k\psi.
\end{align}

It follows that
\begin{align}\label{3.9}
\sum_{1\leq i, k\leq n}c^{kl}u_{ik}\gamma^i=&-\sum_{1\leq i,k\leq n}c^{kl}u_i(\gamma^i)_k+\sum_{1\leq k\leq n}c^{kl}D_k\psi.
\end{align}

Inserting \eqref{3.9} into \eqref{3.7} and combining \eqref{2.6}, \eqref{3.2}, we have
\begin{align}\label{3.10}
\begin{split}
|Dw|^2\log|Dw|^2\frac{\partial\varphi}{\partial\gamma}(x_0)
=&(g'(0)+h'\psi)|Dw|^2\log|Dw|^2-2\sum_{1\leq i,k,l\leq n}c^{kl}u_iu_l(\gamma^i)_k.
\end{split}
\end{align}

From \eqref{3.3}, we obtain
\begin{align}\label{3.11}
|Dw|^2=&|Du|^2-\psi^2\quad\text{on}\quad\partial\Omega.
\end{align}

Assume $|Du|(x_0)\ge \sqrt{100+2|\psi|^2_{C^0(\overline\Omega\times[-M_0, M_0])}}$, otherwise we get the estimates. At $x_0$,   we have
\begin{align}
\frac{1}{2}|Du|^2\leq&|Dw|^2\leq |Du|^2,\label{3.12}\\
|Dw|^2\ge&50.\label{3.13}
\end{align}

Inserting  \eqref{3.12} and \eqref{3.13} into \eqref{3.10}, we have
\begin{align}\label{3.14}
\begin{split}
|Dw|^2\log|Dw|^2\frac{\partial\varphi}{\partial\gamma}(x_0)
\geq&(\alpha_0-|\psi|_{C^0(\overline\Omega\times[-M_0, M_0])}-C_0)|Dw|^2\log|Dw|^2\\
=&|Dw|^2\log|Dw|^2\\
>&0.
\end{split}
\end{align}
On the other hand, from Hopf Lemma, we have
$$\frac{\partial\varphi}{\partial\gamma}(x_0)\leq 0,$$
it is a contradiction to \eqref{3.14}.

Then we have
\begin{align}\label{3.15}
|Du|(x_0)\leq\sqrt{100+2|\psi|^2_{C^0(\overline\Omega\times[-M_0, M_0])}}.
\end{align}

{\bf Case 2.} $ x_0\in \partial\Omega_{\mu_{0}}\bigcap\Omega$. This is due to interior gradient estimates. From Remark~\ref{Rem1.1}, we have
 \begin{align}\label{3.16}
\sup_{\partial\Omega_{\mu_0}\bigcap\Omega}|Du|\leq \tilde{M}_1.
\end{align}
where $\tilde{M}_1$ is a positive constant depending only on $n, M_0, \mu_0, L_1$.\par

{\bf Case 3.} $x_0\in\Omega_{\mu_{0}}$. \par

In this case, $x_0$ is a critical point of $\varphi$. We choose the normal coordinate at $x_0$, by rotating the coordinate system
suitably, we may assume that $u_i(x_0)=0,\,2\leq i\leq n$ and $u_1(x_0)=|Du|>0$. And we can further assume that the matrix $(u_{ij}(x_0))(2\leq i,j\leq n)$ is diagonal. Let $$\mu_2 \le \frac{1}{10
0L_2}$$ such that
 \begin{align}\label{3.16a}
 |\psi_u| \mu_2 \le \frac{1}{100
}, \quad \text{then}\quad \frac{99}{100}\le 1-\psi_u \mu_2 \le
 \frac{101}{100}.
\end{align}

 We can choose $$ \mu_0=\frac{1}{2} \min\{ \mu_1,\mu_2,  1 \}.$$

In order to simplify the calculations, we let
\begin{align*}
w=&u-G,\quad G=\psi(x,u)d.
\end{align*}
Then we have
\begin{align}
w_k=(1-G_u)u_k-G_{x_k}.\label{3wk}
\end{align}
Since at $x_0$,
\begin{align}
|Dw|^2=&w_1^2+\sum_{2\leq i\leq n}w_i^2,\label{3Aa}\\
w_i=&-G_{x_i}=-\psi_{x_i}d-\psi\gamma^i,\quad i=2,\ldots,n,\label{3wi2}\\
w_1=&(1-G_u)u_1-G_{x_1}=(1-G_u)u_1 -\psi_{x_1}d-\psi\gamma^1.\label{3wi2a}
\end{align}
So from the above relation,  at $x_0$, we can assume
\begin{align}
u_1 =  |Du|(x_0)\ge \sqrt{3000(1+ |\psi|^2_{C^1(
\overline\Omega\times[-M_0, M_0])})},\label{3wk}
\end{align}
then
\begin{align}\label{3.18}
\frac{9}{10}u_1^2\leq|Dw|^2\leq \frac{11}{10}u_1^2, \quad \frac{9}{10}u_1^2\leq w_1^2\leq \frac{11}{10}u_1^2
\end{align}
and by
 the choice of $\mu_0$ and  \eqref{3.16a}, we have
\begin{align}
\frac{99}{100}
\le 1-G_u \le\frac {101}{100}.
\label{wi2aa}
\end{align}

From the above choice, we shall prove Theorem~\ref{Thm1.1} with three steps, as we mentioned before, all the calculations will be done at the fixed point $x_0$.

{\bf Step 1:} We first get the formula \eqref{3aijvarphiijb}.\par

Taking the first  derivatives of $\varphi$,
\begin{align}\label{3varphii}
\varphi_i=&\frac{(|Dw|^2)_i}{|Dw|^2\log|Dw|^2}+h'u_i+g'\gamma^i.
\end{align}
From $\varphi_i(x_0)=0$,  we have
\begin{align}\label{3varphii=0}
(|Dw|^2)_i=-|Dw|^2\log|Dw|^2(h'u_i+g'\gamma^i).
\end{align}
Take the derivatives again for  $\varphi_i$,
\begin{align}\label{3varphiija}
\begin{split}
\varphi_{ij}
=&\frac{(|Dw|^2)_{ij}}{|Dw|^2\log|Dw|^2}-(1+\log|Dw|^2)\frac{(|Dw|^2)_i(|Dw|^2)_j}{(|Dw|^2\log|Dw|^2)^2}\\&
+h'u_{ij}
+h''u_iu_j+g''\gamma^i\gamma^j+g'(\gamma^i)_j.
\end{split}
\end{align}
Using \eqref{3varphii=0}, it follows that
\begin{align}\label{3varphiijb}
\begin{split}
\varphi_{ij}
=&\frac{(|Dw|^2)_{ij}}{|Dw|^2\log|Dw|^2}+h'u_{ij}+\big[h''-(1+\log|Dw|^2)h'^2\big]u_iu_j\\&
+\big[g''-(1+\log|Dw|^2)g'^2\big]\gamma^i\gamma^j
-(1+\log|Dw|^2)h'g'(\gamma^iu_j+\gamma^ju_i)+g'(\gamma^i)_j.
\end{split}
\end{align}
Then we get
\begin{align}\label{3aijvarphiija}
\begin{split}
0\geq \sum_{1\leq i,j\leq n}a^{ij}\varphi_{ij}
=:&I_1+I_2,
\end{split}
\end{align}
where
\begin{align}\label{3I1a}
I_1=\frac{1}{|Dw|^2\log|Dw|^2}\sum_{1\leq i,j\leq n}a^{ij}(|Dw|^2)_{ij},
\end{align}
and
\begin{align}
I_2=&\sum_{1\leq i,j\leq n}a^{ij}\bigg\{h' u_{ij}+\big[h''-(1+\log|Dw|^2)h'^2\big]u_iu_j+\big[g''-(1+\log|Dw|^2)g'^2\big]\gamma^i\gamma^j\notag\\
&\qquad\qquad\quad-2(1+\log|Dw|^2)h'g'\gamma^iu_j+g'(\gamma^i)_j\bigg\}.\label{3I2a}
\end{align}
From the choice of the coordinate, we have
\begin{align}\label{3aij}
a^{11}=1,\, a^{ii}=v^2=1+u_1^2 \,\,(2\leq i\leq n ), a^{ij}=0\,\, (i\neq j,1\leq i,j\leq n).
\end{align}

Now we first treat $I_2$.

From the choice of the coordinate and the equations \eqref{2.5}, \eqref{3aij}, we have
\begin{align}
I_2
=&h'f v^3-h'^2u_1^2\log|Dw|^2+(h''-h'^2)u_1^2+\big[g''-(1+\log|Dw|^2)g'^2\big]\sum_{1\leq i\leq n}a^{ii}(\gamma^i)^2\notag\\
&-2(1+\log|Dw|^2)h'g'\gamma^1u_1+g'\sum_{1\leq i\leq n}a^{ii}(\gamma^i)_i\notag\\
=&h'f v^3-(h'^2+c^{11}g'^2) u_1^2\log|Dw|^2+\big[h''-h'^2+c^{11}(g''-g'^2)+g'\sum_{2\leq i\leq n}(\gamma^i)_i\big]u_1^2\notag\\
&-2h'g'\gamma^1 u_1\log|Dw|^2-2h'g'\gamma^1 u_1+g''-g'^2
+g'\sum_{1\leq i\leq n}(\gamma^i)_i.\label{3I2b}
\end{align}
Hence from \eqref{3.1}, we have
\begin{align}
I_2
=&f v^3-(1+c^{11}\alpha_0^2) u_1^2\log|Dw|^2+\big[-1-c^{11}\alpha_0^2+\alpha_0\sum_{2\leq i\leq n}(\gamma^i)_i\big] u_1^2\notag\\
&-2\alpha_0\gamma^1 u_1\log|Dw|^2-2\alpha_0\gamma^1u_1-\alpha_0^2
+\alpha_0\sum_{1\leq i\leq n}(\gamma^i)_i\notag\\
\geq&f v^3-(1+c^{11}\alpha_0^2) u_1^2\log|Dw|^2-C_1u_1^2.\label{3I2c}
\end{align}
here  $C_1$  is a positive constant  depending  only on $n, \Omega, M_0, \mu_0, L_2$.\par

Next, we calculate $I_1$ and get the formula \eqref{3I1c}. \par

Taking the first  derivatives of $|Dw|^2$, we have
\begin{align}\label{3|Dw|^2i}
\begin{split}
(|Dw|^2)_i=&2\sum_{1\leq k\leq n}w_kw_{ki}.
\end{split}
\end{align}
Taking the derivatives of $|Dw|^2$ once more, we have
\begin{align}\label{3|Dw|^2ij}
\begin{split}
(|Dw|^2)_{ij}=&2\sum_{1\leq k\leq n}w_kw_{kij}+2\sum_{1\leq k\leq n}w_{ki}w_{kj}.
\end{split}
\end{align}
By \eqref{3I1a} and \eqref{3|Dw|^2ij}, we can rewrite $I_1$ as
\begin{align}
I_1
=&\frac{1}{|Dw|^2\log|Dw|^2}\big[2\sum_{1\leq i,j,k\leq n}w_ka^{ij}w_{ijk}+2\sum_{1\leq i,j,k\leq n}a^{ij}w_{ki}w_{kj}\big]\notag\\
=:&\frac{1}{|Dw|^2\log|Dw|^2}\big[I_{11}+I_{12}\big].\label{3I1b}
\end{align}
In the following, we shall deal with $I_{11}$ and $I_{12}$ respectively. \par

For the term $I_{11}$. As we have
let
\begin{align}
w=&u-G,\quad G=\psi(x,u)d.\label{3w}
\end{align}
Then we have
\begin{align}
w_k=&(1-G_u)u_k-G_{x_k},\notag\\
w_{ki}=&(1-G_u)u_{ki}-G_{uu}u_ku_i-G_{ux_i}u_k-G_{ux_k}u_i-G_{x_kx_i},\label{3wki}\\
w_{kij}=&(1-G_u)u_{kij}-G_{uu}(u_{ki}u_j+u_{kj}u_i+u_{ij}u_k)\notag\\
&-G_{ux_i}u_{kj}-G_{ux_j}u_{ki}-G_{ux_k}u_{ij}\notag\\
&-G_{uuu}u_ku_iu_j-G_{uux_i}u_ju_k-G_{uux_j}u_iu_k-G_{uux_k}u_iu_j\notag\\
&-G_{ux_ix_j}u_k-G_{ux_kx_j}u_i-G_{ux_ix_k}u_j-G_{x_ix_jx_k}.\label{3wkij}
\end{align}

So from the choice of the coordinate and the equations \eqref{2.5}, \eqref{3aij}, we have
\begin{align}
\sum_{1\leq i,j,k\leq n}w_ka^{ij}w_{ijk}=&\sum_{1\leq i,j,k\leq n}\big[(1-G_u)w_ka^{ij}u_{ijk}-G_{uu}w_ku_ka^{ij}u_{ij}
-2G_{uu}w_ka^{ij}u_{ki}u_j\notag\\
&\qquad\qquad-w_kG_{ux_k}a^{ij}u_{ij}-2w_ka^{ij}G_{ux_i}u_{kj}-G_{uuu}w_ku_ka^{ij}u_iu_j\notag\\
&\qquad\qquad-2w_ku_ka^{ij}G_{uux_i}u_j-w_kG_{uux_k}a^{ij}u_iu_j
-w_ku_ka^{ij}G_{ux_ix_j}\notag\\
&\qquad\qquad-2w_ka^{ij}G_{ux_ix_k}u_j-w_ka^{ij}G_{x_ix_jx_k}\big]\notag\\
=&(1-G_u)\sum_{1\leq i,j,k\leq n}w_ka^{ij}u_{ijk}-2(G_{uu}u_1+G_{ux_1})\sum_{1\leq k\leq n}w_ku_{k1}\notag\\
&-2v^2\sum_{2\leq i\leq n}G_{ux_i}\sum_{1\leq k\leq n}w_ku_{ki}-G_{uu}fv^3u_1w_1
-G_{uuu}u_1^3w_1\notag\\&
-v^2u_1w_1\sum_{2\leq i\leq n}G_{ux_ix_i}-fv^3\sum_{1\leq k\leq n}w_kG_{ux_k}-2G_{uux_1}u_1^2w_1\notag\\&-u_1^2\sum_{1\leq k\leq n}w_kG_{uux_k}
-v^2\sum_{1\leq k\leq n}w_k\sum_{2\leq i\leq n}G_{x_ix_ix_k}-G_{ux_1x_1}u_1w_1\notag\\
&-2u_1\sum_{1\leq k\leq n}G_{ux_1x_k}w_k-\sum_{1\leq k\leq n}G_{x_1x_1x_k}w_k.\label{3wkaijwkij}
\end{align}
By the equation \eqref{2.5}, we have
\begin{align}\label{3u11a}
\begin{split}
u_{11}=&fv^3-v^2\sum_{2\leq i\leq n}u_{ii},
\end{split}
\end{align}
and
\begin{align}\label{3Delta u}
\begin{split}
\Delta u=&fv+\frac{u_1^2}{v^2}u_{11}.
\end{split}
\end{align}
Differentiating \eqref{2.5}, we have
\begin{align}\label{3aijuijka}
\begin{split}
\sum_{1\leq i,j\leq n}a^{ij}u_{ijk}=&-\sum_{1\leq i,j,l\leq n}a^{ij}_{p_l}u_{lk}u_{ij}+v^3D_{k}f+3fv^2v_k.
\end{split}
\end{align}
From \eqref{2.4}, we have
\begin{align}\label{3aijpl}
\begin{split}
a^{ij}_{p_l}=&2u_l\delta_{ij}-\delta_{il}u_j-\delta_{jl}u_i.
\end{split}
\end{align}
By the definition of $v$, we have
\begin{align}\label{3vvk}
\begin{split}
v v_k=&u_1u_{1k}.
\end{split}
\end{align}
Since
\begin{align}\label{3Dkf}
\begin{split}
D_kf=&f_uu_k+f_{x_k},
\end{split}
\end{align}
 from \eqref{3aijpl}, \eqref{3u11a} and \eqref{3vvk}, we have
\begin{align}\label{3aijuijkb}
\begin{split}
\sum_{1\leq i,j\leq n}a^{ij}u_{ijk}
=&-2u_1u_{1k}\Delta u+2u_1\sum_{1\leq i\leq n}u_{1i}u_{ik}+v^3D_{k}f+3fvu_1u_{1k},\\
=&\frac{2u_1}{v^2}u_{11}u_{1k}+2u_1\sum_{2\leq i\leq n}u_{1i}u_{ik}+v^3D_{k}f+fvu_1u_{1k}\\
=&\frac{2u_1}{v^2}u_{11}u_{1k}+2u_1\sum_{2\leq i\leq n}u_{1i}u_{ik}+f_uv^3u_k+f_{x_k}v^3+fvu_1u_{1k}.
\end{split}
\end{align}
Inserting \eqref{3aijuijkb} into \eqref{3wkaijwkij}, we rewrite $I_{11}$ as
\begin{align}\label{3I11a}
I_{11}
=&4(1-G_u)\frac{u_1}{v^2} u_{11}\sum_{1\leq k\leq n}w_ku_{k1}+4(1-G_u)u_1\sum_{2\leq i\leq n}u_{1i}\sum_{1\leq k\leq n}w_ku_{ki}\notag\\
&+[2(1-G_u)fu_1v-4G_{uu}u_1-4G_{ux_1}]\sum_{1\leq k\leq n}w_ku_{k1}-4v^2\sum_{2\leq i\leq n}G_{ux_i}\sum_{1\leq k\leq n}w_ku_{ki}\notag\\
&+2(1-G_u)f_uv^3u_1w_1-2G_{uu}fv^3u_1w_1+2(1-G_u)v^3\sum_{1\leq k\leq n}f_{x_k}w_k\notag\\&-2fv^3\sum_{1\leq k\leq n}G_{ux_k}w_k
-2G_{uuu}u_1^3w_1-2v^2u_1w_1\sum_{2\leq i\leq n}G_{ux_ix_i}-4G_{uux_1}u_1^2w_1\notag\\&-2u_1^2\sum_{1\leq k\leq n}w_kG_{uux_k}
-2v^2\sum_{1\leq k\leq n}w_k\sum_{2\leq i\leq n}G_{x_ix_ix_k}-2G_{ux_1x_1}u_1w_1\notag\\
&-4u_1\sum_{1\leq k\leq n}G_{ux_1x_k}w_k-2\sum_{1\leq k\leq n}G_{x_1x_1x_k}w_k.
\end{align}

For the term $I_{12}$: applying \eqref{3aij} and \eqref{3wki}, we have
\begin{align*}
I_{12}
=&2\sum_{1\leq i,k\leq n}a^{ii}w^2_{ki}\notag\\
=&2\sum_{1\leq k\leq n}w^2_{k1}+2v^2\sum_{1\leq k\leq n}\sum_{2\leq i\leq n}w^2_{ki}\notag\\
=&2w_{11}^2+2(v^2+1)\sum_{2\leq i\leq n}w^2_{1i}+2v^2\sum_{2\leq i,k\leq n}w^2_{ki}\notag\\
=&2[(1-G_u)u_{11}-(G_{uu}u_1^2+2G_{ux_1}u_1+G_{x_1x_1})]^2\notag\\
&+2(v^2+1)\sum_{2\leq i\leq n}[(1-G_u)u_{1i}-(G_{ux_i}u_1+G_{x_1x_i})]^2\notag\\
&+2v^2\sum_{2\leq i,k\leq n}[(1-G_u)u_{ki}-G_{x_kx_i}]^2\notag\\
=&2(1-G_u)^2u_{11}^2+2(1-G_u)^2(v^2+1)\sum_{2\leq i\leq n}u_{1i}^2+2(1-G_u)^2v^2\sum_{2\leq i\leq n}u_{ii}^2\notag\\
\end{align*}
\begin{align}
&-4(1-G_u)(G_{uu}u_1^2+2G_{ux_1}u_1+G_{x_1x_1})u_{11}\notag\\
&-4(1-G_u)(v^2+1)\sum_{2\leq i\leq n}(G_{ux_i}u_1+G_{x_1x_i})u_{1i}\notag\\
&-4(1-G_u)v^2\sum_{2\leq i\leq n}G_{x_ix_i}u_{ii}+2(G_{uu}u_1^2+2G_{ux_1}u_1+G_{x_1x_1})^2\notag\\
&+2(v^2+1)\sum_{2\leq i\leq n}(G_{ux_i}u_1+G_{x_1x_i})^2
+2v^2\sum_{2\leq i,k\leq n}G_{x_kx_i}.\label{3I12a}
\end{align}
Combining  \eqref{3I11a}, \eqref{3I12a}, it follows that
\begin{align}
I_1=&\frac{1}{|Dw|^2\log|Dw|^2}\bigg[2(1-G_u)^2u_{11}^2+2(1-G_u)^2(v^2+1)\sum_{2\leq i\leq n}u_{1i}^2\notag\\&+2(1-G_u)^2v^2\sum_{2\leq i\leq n}u_{ii}^2
+4(1-G_u)\frac{u_1}{v^2} u_{11}\sum_{1\leq k\leq n}w_ku_{k1}\notag\\&+4(1-G_u)u_1\sum_{2\leq i\leq n}u_{1i}\sum_{1\leq k\leq n}w_ku_{ki}\notag\\&
+[2(1-G_u)fu_1v-4G_{uu}u_1-4G_{ux_1}]\sum_{1\leq k\leq n}w_ku_{k1}\notag\\&-4v^2\sum_{2\leq i\leq n}G_{ux_i}\sum_{1\leq k\leq n}w_ku_{ki}\notag\\
&-4(1-G_u)(G_{uu}u_1^2+2G_{ux_1}u_1+G_{x_1x_1})u_{11}\notag\\
&-4(1-G_u)(v^2+1)\sum_{2\leq i\leq n}(G_{ux_i}u_1+G_{x_1x_i})u_{1i}\notag\\
&-4(1-G_u)v^2\sum_{2\leq i\leq n}G_{x_ix_i}u_{ii}+2(1-G_u)f_uv^3u_1w_1-2G_{uu}fv^3u_1w_1\notag\\
&+2(G_{uu}u_1^2+2G_{ux_1}u_1+G_{x_1x_1})^2+2(v^2+1)\sum_{2\leq i\leq n}(G_{ux_i}u_1+G_{x_1x_i})^2\notag\\
&+2v^2\sum_{2\leq i,k\leq n}G_{x_kx_i}+2(1-G_u)v^3\sum_{1\leq k\leq n}f_{x_k}w_k-2fv^3\sum_{1\leq k\leq n}G_{ux_k}w_k\notag\\&
-2G_{uuu}u_1^3w_1-2v^2u_1w_1\sum_{2\leq i\leq n}G_{ux_ix_i}-4G_{uux_1}u_1^2w_1-2u_1^2\sum_{1\leq k\leq n}w_kG_{uux_k}\notag\\&
-2v^2\sum_{1\leq k\leq n}w_k\sum_{2\leq i\leq n}G_{x_ix_ix_k}-2G_{ux_1x_1}u_1w_1\notag\\
&-4u_1\sum_{1\leq k\leq n}G_{ux_1x_k}w_k-2\sum_{1\leq k\leq n}G_{x_1x_1x_k}w_k\bigg].\label{3I1c}
\end{align}
Inserting \eqref{3I1c} and \eqref{3I2b} into \eqref{3aijvarphiija}, we  can obtain the following formula
\begin{align}\label{3aijvarphiijb}
\begin{split}
0\geq \sum_{1\leq i,j\leq n}a^{ij}\varphi_{ij}=: Q_1+Q_2+Q_3,
\end{split}
\end{align}
where $Q_1$ contains all the quadratic terms of $u_{ij}$; $Q_2$  is the term which contains all linear terms of $u_{ij}$;  and   the remaining terms are denoted by $Q_3$.
Then we have
\begin{align}\label{3Q1a}
Q_1=&\frac{1}{|Dw|^2\log|Dw|^2}\bigg[ 2(1-G_u)^2u_{11}^2+2(1-G_u)^2(v^2+1)\sum_{2\leq i\leq n}u_{1i}^2\notag\\&
+4(1-G_u)\frac{u_1}{v^2} u_{11}\sum_{1\leq k\leq n}w_ku_{k1}
+4(1-G_u)u_1\sum_{2\leq i\leq n}u_{1i}\sum_{1\leq k\leq n}w_ku_{ki}\notag\\&
+2(1-G_u)^2v^2\sum_{2\leq i\leq n}u_{ii}^2\bigg];
\end{align}
The linear terms of $u_{ij}$ are
\begin{align}\label{3Q2a}
Q_2=&\frac{1}{|Dw|^2\log|Dw|^2}\bigg[\big[2(1-G_u)fu_1v-4G_{uu}u_1-4G_{ux_1}\big]\sum_{1\leq k\leq n}w_ku_{k1}\notag\\&
-4v^2\sum_{2\leq i\leq n}G_{ux_i}\sum_{1\leq k\leq n}w_ku_{ki}-4(1-G_u)(G_{uu}u_1^2+2G_{ux_1}u_1+G_{x_1x_1})u_{11}\notag\\
&-4(1-G_u)(v^2+1)\sum_{2\leq i\leq n}(G_{ux_i}u_1+G_{x_1x_i})u_{1i}-4(1-G_u)v^2\sum_{2\leq i\leq n}G_{x_ix_i}u_{ii}\bigg];
\end{align}
and  the
remaining terms are
\begin{align}
Q_3=&I_2+\frac{1}{|Dw|^2\log|Dw|^2}\bigg[2(1-G_u)f_uv^3u_1w_1-2G_{uu}fv^3u_1w_1\notag\\
&+2(G_{uu}u_1^2+2G_{ux_1}u_1+G_{x_1x_1})^2+2(v^2+1)\sum_{2\leq i\leq n}(G_{ux_i}u_1+G_{x_1x_i})^2\notag\\
&+2v^2\sum_{2\leq i,k\leq n}G_{x_kx_i}+2(1-G_u)v^3\sum_{1\leq k\leq n}f_{x_k}w_k-2fv^3\sum_{1\leq k\leq n}G_{ux_k}w_k\notag\\&
-2G_{uuu}u_1^3w_1-2v^2u_1w_1\sum_{2\leq i\leq n}G_{ux_ix_i}-4G_{uux_1}u_1^2w_1-2u_1^2\sum_{1\leq k\leq n}w_kG_{uux_k}\notag\\&
-2v^2\sum_{1\leq k\leq n}w_k\sum_{2\leq i\leq n}G_{x_ix_ix_k}-2G_{ux_1x_1}u_1w_1\notag\\
&-4u_1\sum_{1\leq k\leq n}G_{ux_1x_k}w_k-2\sum_{1\leq k\leq n}G_{x_1x_1x_k}w_k\bigg].\label{3Q3a}
\end{align}
From the estimate on $I_2$ in \eqref{3I2c}, we have
\begin{align}\label{3Q3b}
Q_3\geq&f v^3-\frac{2fG_{uu}}{|Dw|^2\log|Dw|^2}v^3u_1w_1-(1+c^{11}\alpha_0^2) u_1^2\log|Dw|^2-C_2u_1^2,
\end{align}
in the computation of $Q_3$, we use the relation $f_u\geq 0$, where $C_2$ is a positive constant which depends only on $n, \Omega, M_0, \mu_0, L_1, L_2$.

{\bf Step 2:} In this step we shall treat the terms $Q_1, Q_2$ , using the first order derivative condition
 $$\varphi_i(x_0)=0,$$
  and let
  \begin{align}\label{3A}
A=|Dw|^2\log|Dw|^2.
\end{align}

By \eqref{3varphii=0}  and  \eqref{3|Dw|^2i},  we have
\begin{align}\label{3wkwki=1}
\begin{split}
\sum_{1\leq k\leq n}w_kw_{ki}=&-\frac{h'}{2}|Dw|^2\log|Dw|^2u_i-\frac{g'\gamma^i}{2}|Dw|^2\log|Dw|^2\\
=&-\frac{h'}{2}Au_i-\frac{g'\gamma^i}{2}A,\qquad i=1,2,\ldots,n.
\end{split}
\end{align}

Putting \eqref{3wki} into \eqref{3wkwki=1},  we get
\begin{align}\label{3wkuki=1}
\begin{split}
(1-G_u)\sum_{1\leq k\leq n} w_ku_{ki}
=&-\frac{h'}{2}Au_i-\frac{g'\gamma^i}{2}A+G_{uu}\sum_{1\leq k\leq n} w_ku_ku_i+\sum_{1\leq k\leq n} w_ku_kG_{ux_i}\\
&+\sum_{1\leq k\leq n} w_kG_{ux_k}u_i+\sum_{1\leq k\leq n} w_kG_{x_kx_i}\\
=&-\frac{h'}{2}Au_i-\frac{g'\gamma^i}{2}A+G_{uu}w_1u_1u_i+ w_1u_1G_{ux_i}\\
&+\sum_{1\leq k\leq n} w_kG_{ux_k}u_i+\sum_{1\leq k\leq n} w_kG_{x_kx_i},\quad i=1,2,\ldots,n.
\end{split}
\end{align}

By \eqref{3wkuki=1}, we have
\begin{align}\label{3wkuk1}
\begin{split}
(1-G_u)\sum_{1\leq k\leq n} w_ku_{k1}
=&-\frac{h'}{2}Au_1-\frac{g'\gamma^1}{2}A+G_{uu}u_1^2w_1+G_{ux_1}w_1u_1\\
&+\sum_{1\leq k\leq n} w_kG_{ux_k}u_1+\sum_{1\leq k\leq n} w_kG_{x_kx_1},
\end{split}
\end{align}
and
\begin{align}\label{3wkuki=2}
\begin{split}
(1-G_u)\sum_{1\leq k\leq n} w_ku_{ki}
=&-\frac{g'\gamma^i}{2}A+G_{ux_i}w_1u_1+\sum_{1\leq k\leq n} w_kG_{x_kx_i},\quad i=2,\ldots,n.
\end{split}
\end{align}

Through \eqref{3wkuki=2}  and the choice of the coordinate at $x_0$, we have
\begin{align}\label{3u1i}
\begin{split}
u_{1i}
=&-\frac{1}{w_1}w_iu_{ii}-\frac{g'\gamma^i}{2(1-G_u)}\frac{A}{w_1}+\frac{G_{ux_i}}{1-G_u}u_1+\frac{1}{(1-G_u)w_1}\sum_{1\leq k\leq n} w_kG_{x_kx_i},\\
&\hspace{200pt} i=2,3,\ldots,n.
\end{split}
\end{align}
Using \eqref{3wkuk1}  and \eqref{3u1i},  it follows that
\begin{align}\label{3u11b}
\begin{split}
u_{11}
=&\sum_{2\leq i\leq n}\frac{w_i^2}{w_1^2}u_{ii}-\frac{h'}{2(1-G_u)}\frac{Au_1}{w_1}+\frac{G_{uu}}{1-G_u}u_1^2
-\frac{g'\gamma^1}{2(1-G_u)}\frac{A}{w_1}\\
&+\frac{G_{ux_1}}{1-G_u}u_1+\frac{u_1}{(1-G_u)w_1}\sum_{1\leq k\leq n} w_kG_{ux_k}+\frac{g'}{2(1-G_u)}\sum_{2\leq i\leq n}w_i\gamma^i\frac{A}{w_1^2}\\
&-\frac{u_1}{(1-G_u)w_1}\sum_{2\leq i\leq n}G_{ux_i}w_i+\frac{1}{(1-G_u)w_1}\sum_{1\leq k\leq n} w_kG_{x_kx_1}\\&
-\frac{1}{(1-G_u)w_1^2}\sum_{1\leq k,i\leq n} w_iw_kG_{x_kx_i}\\
=&:\sum_{2\leq i\leq n}\frac{w_i^2}{w_1^2}u_{ii}-\frac{h'}{2(1-G_u)}\frac{Au_1}{w_1}+\frac{D}{1-G_u},
\end{split}
\end{align}
 where we have let
 \begin{align}\label{3Da}
\begin{split}
D=&G_{uu}u_1^2
-\frac{g'\gamma^1}{2}\frac{A}{w_1}+G_{ux_1}u_1+\frac{u_1}{w_1}\sum_{1\leq k\leq n} w_kG_{ux_k}+\frac{g'}{2}\sum_{2\leq i\leq n}w_i\gamma^i\frac{A}{w_1^2}\\
&-\frac{u_1}{w_1}\sum_{2\leq i\leq n}G_{ux_i}w_i+\frac{1}{w_1}\sum_{1\leq k\leq n} w_kG_{x_kx_1}
-\frac{1}{w_1^2}\sum_{1\leq k,i\leq n} w_iw_kG_{x_kx_i}.
\end{split}
\end{align}
It follows that
\begin{align}\label{3Db}
\begin{split}
|D|
\le C_3 u_1^2.
\end{split}
\end{align}

By \eqref{3u11a}  and \eqref{3u11b},  we have
\begin{align}\label{3uii=2}
\begin{split}
\sum_{2\leq i\leq n}(v^2+\frac{w_i^2}{w_1^2})u_{ii}
=f v^3+\frac{h'}{2(1-G_u)}\frac{Au_1}{w_1}-\frac{D}{1-G_u}.
\end{split}
\end{align}

Now we use the formulas \eqref{3wkuk1}-\eqref{3u11b}  to treat  each term in $Q_1, Q_2$.  At first, we treat  the first four terms of $Q_1$ in \eqref{3Q1a}, and get \eqref{3u11^2}-\eqref{3u1iwkuki}.

By \eqref{3u11b}, we have
\begin{align}\label{3u11^2}
2(1-G_u)^2u_{11}^2=&2(1-G_u)^2\big[\sum_{2\leq i\leq n}\frac{w_i^2}{w_1^2}u_{ii}-\frac{h'}{2(1-G_u)}\frac{Au_1}{w_1}+\frac{D}{1-G_u}\big]^2\notag\\
=&2(1-G_u)^2(\sum_{2\leq i\leq n}\frac{w_i^2}{w_1^2}u_{ii})^2-2(1-G_u)(h'\frac{Au_1}{w_1^3}-2\frac{D}{w_1^2})\sum_{2\leq i\leq n}w_i^2u_{ii}\notag\\
&+\frac{h'^2}{2}\frac{A^2u_1^2}{w_1^2}-2h'\frac{Au_1}{w_1}D+2D^2
\notag\\
=&2(1-G_u)^2(\sum_{2\leq i\leq n}\frac{w_i^2}{w_1^2}u_{ii})^2-2(1-G_u)(h'\frac{Au_1}{w_1^3}-2\frac{D}{w_1^2})\sum_{2\leq i\leq n}w_i^2u_{ii}\notag\\
&+\frac{h'^2}{2}\frac{A^2u_1^2}{w_1^2}+AO(u_1^2).
\end{align}

By \eqref{3u1i}, we have
\begin{align}\label{3u1i^2}
&2(1-G_u)^2(v^2+1)\sum_{2\leq i\leq n}u_{1i}^2\notag\\
=&2(v^2+1)\sum_{2\leq i\leq n}\big[(1-G_u)\frac{w_i}{w_1}u_{ii}+\frac{g'\gamma^i}{2}\frac{A}{w_1}-G_{ux_i}u_1-\frac{1}{w_1}\sum_{1\leq k\leq n} w_kG_{x_kx_i}\big]^2\notag\\
=&2(1-G_u)^2(v^2+1)\sum_{2\leq i\leq n}\frac{w_i^2}{w_1^2}u_{ii}^2+2(1-G_u)g'A\frac{v^2+1}{w_1^2}\sum_{2\leq i\leq n}w_i\gamma^iu_{ii}\notag\\
&-4(1-G_u)\frac{v^2+1}{w_1^2}\sum_{2\leq i\leq n}(G_{ux_i}u_1w_1+\sum_{1\leq k\leq n} w_kG_{x_kx_i})w_iu_{ii}\notag\\
&+\frac{c^{11}}{2}g'^2\frac{v^2+1}{w_1^2}A^2-2g'A\frac{v^2+1}{w_1^2}\sum_{2\leq i\leq n}(G_{ux_i}u_1w_1+\sum_{1\leq k\leq n} w_kG_{x_kx_i})\gamma^{i}\notag\\
&+\frac{2(v^2+1)}{w_1^2}\sum_{2\leq i\leq n}(G_{ux_i}u_1w_1+\sum_{1\leq k\leq n} w_kG_{x_kx_i})^2
\notag\\
=&2(1-G_u)^2(v^2+1)\sum_{2\leq i\leq n}\frac{w_i^2}{w_1^2}u_{ii}^2+2(1-G_u)g'A\frac{v^2+1}{w_1^2}\sum_{2\leq i\leq n}w_i\gamma^iu_{ii}\notag\\
&-4(1-G_u)\frac{v^2+1}{w_1^2}\sum_{2\leq i\leq n}(G_{ux_i}u_1w_1+\sum_{1\leq k\leq n} w_kG_{x_kx_i})w_iu_{ii}\notag\\
&+\frac{c^{11}}{2}g'^2\frac{v^2+1}{w_1^2}A^2+AO(u_1^2).
\end{align}

By \eqref{3u11b} and \eqref{3wkuk1}, we have
\begin{align}\label{3u11wkuk1}
&4(1-G_u)\frac{u_1}{v^2} u_{11}\sum_{1\leq k\leq n}w_ku_{k1}\notag\\
=&\frac{4u_1}{v^2}\big[\sum_{2\leq i\leq n}\frac{w_i^2}{w_1^2}u_{ii}-\frac{h'}{2(1-G_u)}\frac{Au_1}{w_1}+\frac{D}{1-G_u}\big]
\big[-\frac{h'}{2}Au_1-\frac{g'\gamma^1}{2}A+G_{uu}u_1^2w_1\notag\\&+G_{ux_1}w_1u_1
+\sum_{1\leq k\leq n} w_kG_{ux_k}u_1+\sum_{1\leq k\leq n} w_kG_{x_kx_1}\big]\notag\\
=&\big[-2h'\frac{Au_1^2}{v^2w_1^2}+4G_{uu}\frac{u_1^3}{v^2w_1}-2g'\gamma^1\frac{Au_1}{v^2w_1^2}+4G_{ux_1}\frac{u_1^2}{v^2w_1}
+4\sum_{1\leq k\leq n} w_kG_{ux_k}\frac{u_1^2}{v^2w_1^2}\notag\\&
+4\sum_{1\leq k\leq n} w_kG_{x_kx_1}\frac{u_1}{v^2w_1^2}\big]\sum_{2\leq i\leq n}w_i^2u_{ii}+\frac{h'^2}{1-G_u}\frac{A^2u_1^3}{v^2w_1}-\frac{2h'}{1-G_u}\frac{Au_1^2}{v^2}D\notag\\&
-\frac{4u_1}{v^2}[\frac{h'}{2(1-G_u)}\frac{Au_1}{w_1}-\frac{D}{1-G_u}][\frac{g'\gamma^1}{2}A-G_{uu}u_1^2w_1-G_{ux_1}w_1u_1\notag\\&
-\sum_{1\leq k\leq n} w_kG_{ux_k}u_1-\sum_{1\leq k\leq n} w_kG_{x_kx_1}]
\notag\\
=&\big[-2h'\frac{Au_1^2}{v^2w_1^2}+4G_{uu}\frac{u_1^3}{v^2w_1}-2g'\gamma^1\frac{Au_1}{v^2w_1^2}+4G_{ux_1}\frac{u_1^2}{v^2w_1}
+4\sum_{1\leq k\leq n} w_kG_{ux_k}\frac{u_1^2}{v^2w_1^2}\notag\\&
+4\sum_{1\leq k\leq n} w_kG_{x_kx_1}\frac{u_1}{v^2w_1^2}\big]\sum_{2\leq i\leq n}w_i^2u_{ii}+\frac{h'^2}{1-G_u}\frac{A^2u_1^3}{v^2w_1}+AO(u_1^2).
\end{align}

By \eqref{3u1i} and \eqref{3wkuki=2}, we have
\begin{align*}
&4(1-G_u)u_1\sum_{2\leq i\leq n}u_{1i}\sum_{1\leq k\leq n}w_ku_{ki}\notag\\
=&4u_1\sum_{2\leq i\leq n}\big[-\frac{g'\gamma^i}{2}A+G_{ux_i}w_1u_1+\sum_{1\leq k\leq n} w_kG_{x_kx_i}\big]\big[-\frac{1}{w_1}w_iu_{ii}-\frac{g'\gamma^i}{2(1-G_u)}\frac{A}{w_1}\notag\\&+\frac{G_{ux_i}}{1-G_u}u_1+\frac{1}{(1-G_u)w_1}\sum_{1\leq k\leq n} w_kG_{x_kx_i}\big]\notag\\
=&2g'\frac{Au_1}{w_1}\sum_{2\leq i\leq n}w_i\gamma^iu_{ii}-4u_1^2\sum_{2\leq i\leq n}w_iG_{ux_i}u_{ii}
-\frac{4u_1}{w_1}\sum_{1\leq k\leq n}w_k\sum_{2\leq i\leq n}G_{x_kx_i}w_iu_{ii}\notag\\
&
+\frac{c^{11}g'^2}{1-G_u}\frac{A^2u_1}{w_1}-\sum_{2\leq i\leq n}
\big[2g'\gamma^iAu_1-4G_{ux_i}w_1u_1^2-4u_1\sum_{1\leq k\leq n} w_kG_{x_kx_i}\big]\big[\frac{G_{ux_i}}{1-G_u}u_1\notag\\&+\frac{1}{(1-G_u)w_1}\sum_{1\leq k\leq n} w_kG_{x_kx_i}\big]-\frac{2g'}{1-G_u}\frac{Au_1}{w_1}\sum_{2\leq i\leq n}\big[G_{ux_i}u_1w_1+\sum_{1\leq k\leq n} w_kG_{x_kx_i}\big]\gamma^i
\notag\\
\end{align*}
\begin{align}
=&2g'\frac{Au_1}{w_1}\sum_{2\leq i\leq n}w_i\gamma^iu_{ii}-4u_1^2\sum_{2\leq i\leq n}w_iG_{ux_i}u_{ii}
-\frac{4u_1}{w_1}\sum_{1\leq k\leq n}w_k\sum_{2\leq i\leq n}G_{x_kx_i}w_iu_{ii}\notag\\&
+\frac{c^{11}g'^2}{1-G_u}\frac{A^2u_1}{w_1}+AO(u_1^2).\label{3u1iwkuki}
\end{align}

Now we  treat  the first four terms of $Q_2$ in \eqref{3Q2a}, and get \eqref{3wkuk1a}-\eqref{3u1ia}.

From \eqref{3wkuk1}, we get
\begin{align}\label{3wkuk1a}
&\big[2(1-G_u)fu_1v-4G_{uu}u_1-4G_{ux_1}\big]\sum_{1\leq k\leq n}w_ku_{k1}\notag\\
=&\frac{1}{1-G_u}\big[2(1-G_u)fu_1v-4G_{uu}u_1-4G_{ux_1}\big]\big[-\frac{h'}{2}Au_1+G_{uu}u_1^2w_1\notag\\&-\frac{g'\gamma^1}{2}A+G_{ux_1}w_1u_1
+\sum_{1\leq k\leq n} w_kG_{ux_k}u_1+\sum_{1\leq k\leq n} w_kG_{x_kx_1}\big]\notag\\
=&-h'fAvu_1^2+2G_{uu}fvu_1^3w_1-fg'\gamma^1Avu_1\notag\\
&+2fvu_1[G_{ux_1}w_1u_1
+\sum_{1\leq k\leq n} w_kG_{ux_k}u_1+\sum_{1\leq k\leq n} w_kG_{x_kx_1}]\notag\\
&-\frac{4}{1-G_u}(G_{uu}u_1+G_{ux_1})\big[-\frac{h'}{2}Au_1+G_{uu}u_1^2w_1\notag\\&-\frac{g'\gamma^1}{2}A+G_{ux_1}w_1u_1
+\sum_{1\leq k\leq n} w_kG_{ux_k}u_1+\sum_{1\leq k\leq n} w_kG_{x_kx_1}\big]
\notag\\
=&-h'fvu_1^2+2G_{uu}fvu_1^3w_1+AO(u_1^2).
\end{align}

From \eqref{3wkuki=2}, we have
\begin{align}\label{3wkuki=2a}
&-4v^2\sum_{2\leq i\leq n}G_{ux_i}\sum_{1\leq k\leq n}w_ku_{ki}\notag\\
=&-\frac{4v^2}{1-G_u}\sum_{2\leq i\leq n}G_{ux_i}\big[-\frac{g'\gamma^i}{2}A+G_{ux_i}w_1u_1+\sum_{1\leq k\leq n} w_kG_{x_kx_i}\big]\notag\\
=&\frac{2g'}{1-G_u}\sum_{2\leq i\leq n}G_{ux_i}\gamma^iAv^2-\frac{4}{1-G_u}\sum_{2\leq i\leq n}G^2_{ux_i}v^2u_1w_1\notag\\
&-\frac{4}{1-G_u}\sum_{2\leq i\leq n}G_{ux_i}\sum_{1\leq k\leq n} w_kG_{x_kx_i}v^2
\notag\\
=&AO(u_1^2).
\end{align}

From \eqref{3u11b}, we have
\begin{align}\label{3u11c}
&-4(1-G_u)(G_{uu}u_1^2+2G_{ux_1}u_1+G_{x_1x_1})u_{11}\notag\\
=&-4(G_{uu}u_1^2+2G_{ux_1}u_1+G_{x_1x_1})\big[(1-G_u)\sum_{2\leq i\leq n}\frac{w_i^2}{w_1^2}u_{ii}-\frac{h'}{2}\frac{Au_1}{w_1}+D\big]\notag\\
=&-4(1-G_u)(G_{uu}u_1^2+2G_{ux_1}u_1+G_{x_1x_1})\sum_{2\leq i\leq n}\frac{w_i^2}{w_1^2}u_{ii}\notag\\
&+4(G_{uu}u_1^2+2G_{ux_1}u_1+G_{x_1x_1})(\frac{h'}{2}\frac{Au_1}{w_1}-D)
\notag\\
=&-4(1-G_u)(G_{uu}u_1^2+2G_{ux_1}u_1+G_{x_1x_1})\sum_{2\leq i\leq n}\frac{w_i^2}{w_1^2}u_{ii}+AO(u_1^2).
\end{align}

From \eqref{3u1i}, we have
\begin{align}\label{3u1ia}
&-4(1-G_u)(v^2+1)\sum_{2\leq i\leq n}(G_{ux_i}u_1+G_{x_1x_i})u_{1i}\notag\\
=&4(v^2+1)\sum_{2\leq i\leq n}(G_{ux_i}u_1+G_{x_1x_i})\big[\frac{1-G_u}{w_1}w_iu_{ii}+\frac{g'\gamma^i}{2}\frac{A}{w_1}\notag\\&
-G_{ux_i}u_1-\frac{1}{w_1}\sum_{1\leq k\leq n} w_kG_{x_kx_i}\big]\notag\\
=&4(1-G_u)\frac{v^2+1}{w_1}\sum_{2\leq i\leq n}(G_{ux_i}u_1+G_{x_1x_i})w_iu_{ii}\notag\\
&+4(v^2+1)\sum_{2\leq i\leq n}(G_{ux_i}u_1+G_{x_1x_i})\big[\frac{g'\gamma^i}{2}\frac{A}{w_1}
-G_{ux_i}u_1-\frac{1}{w_1}\sum_{1\leq k\leq n} w_kG_{x_kx_i}\big]
\notag\\
=&4(1-G_u)\frac{v^2+1}{w_1}\sum_{2\leq i\leq n}(G_{ux_i}u_1+G_{x_1x_i})w_iu_{ii}+AO(u_1^2).
\end{align}

We treat the term $Q_1$ using the relations \eqref{3u11^2}-\eqref{3u1iwkuki} and  use the formulas \eqref{3wkuk1a}-\eqref{3u1ia}  to treat the term $Q_2$.
 By the formula on $Q_3$ in \eqref{3Q3a},  we can get the following new formula of \eqref{3aijvarphiijb},
\begin{align}\label{3aijvarphiijc}
0\geq \sum_{1\leq i,j\leq n}a^{ij}\varphi_{ij}=: J_1+J_2,
\end{align}
where $J_1$  only contains the terms with $u_{ii}$ , the other terms belong to $J_2$. We can write
\begin{align}\label{3J1a}
J_1=:\frac{1}{A}\big[J_{11}+J_{12}\big],
\end{align}
here $J_{11}$  contains the quadratic terms of $u_{ii}\,(i\geq2)$, and  $J_{12}$ is the term including linear terms of $u_{ii}\,(i\geq2)$.
It follows that
\begin{align}\label{3J11a}
J_{11}=&2(1-G_u)^2\big\{(\sum_{2\leq i\leq n}\frac{w_i^2}{w_1^2}u_{ii})^2+(v^2+1)\sum_{2\leq i\leq n}\frac{w_i^2}{w_1^2}u_{ii}^2+v^2\sum_{2\leq i\leq n}u_{ii}^2\big\}\notag\\
=&2(1-G_u)^2\big\{\sum_{2\leq i\leq n}\frac{w_i^4}{w_1^4}u_{ii}^2+2\sum_{2\leq i< j\leq n}\frac{w_i^2w_j^2}{w_1^4}u_{ii}u_{jj}
+\sum_{2\leq i\leq n}[v^2+(v^2+1)\frac{w_i^2}{w_1^2}]u_{ii}^2\big\}\notag\\
=&2(1-G_u)^2\big\{\sum_{2\leq i\leq n}(v^2+\frac{w_i^2}{w_1^2})(1+\frac{w_i^2}{w_1^2})u_{ii}^2+2\sum_{2\leq i< j\leq n}\frac{w_i^2w_j^2}{w_1^4}u_{ii}u_{jj}
\big\}\notag\\
=&2(1-G_u)^2\big\{\sum_{2\leq i\leq n}d_ie_iu_{ii}^2+2\sum_{2\leq i< j\leq n}\frac{w_i^2w_j^2}{w_1^4}u_{ii}u_{jj}
\big\},
\end{align}
where
\begin{align}
d_i=&v^2+\frac{w_i^2}{w_1^2}, \quad i=2,3,\ldots,n,\label{3di}\\
e_i=&1+\frac{w_i^2}{w_1^2},\quad i=2,3,\ldots,n.\label{3ei}
\end{align}
And
\begin{align}\label{3J12a}
J_{12}=&-2(1-G_u)(h'\frac{Au_1}{w_1^3}-2\frac{D}{w_1^2})\sum_{2\leq i\leq n}w_i^2u_{ii}
+2(1-G_u)g'A\frac{v^2+1}{w_1^2}\sum_{2\leq i\leq n}w_i\gamma^iu_{ii}\notag\\
&-4(1-G_u)\frac{v^2+1}{w_1^2}\sum_{2\leq i\leq n}(\sum_{1\leq k\leq n} w_kG_{x_kx_i}-G_{x_1x_i}w_1)w_iu_{ii}\notag\\
&-\big[2h'\frac{Au_1^2}{v^2w_1^2}-4G_{uu}\frac{u_1^3}{v^2w_1}+2g'\gamma^1\frac{Au_1}{v^2w_1^2}-4G_{ux_1}\frac{u_1^2}{v^2w_1}
-4\sum_{1\leq k\leq n} w_kG_{ux_k}\frac{u_1^2}{v^2w_1^2}\notag\\&
-4\sum_{1\leq k\leq n} w_kG_{x_kx_1}\frac{u_1}{v^2w_1^2}\big]\sum_{2\leq i\leq n}w_i^2u_{ii}+2g'\frac{Au_1}{w_1}\sum_{2\leq i\leq n}w_i\gamma^iu_{ii}
-4u_1^2\sum_{2\leq i\leq n}w_iG_{ux_i}u_{ii}\notag\\&
-\frac{4u_1}{w_1}\sum_{1\leq k\leq n}w_k\sum_{2\leq i\leq n}G_{x_kx_i}w_iu_{ii}
-4(1-G_u)v^2\sum_{2\leq i\leq n}G_{x_ix_i}u_{ii}\notag\\&-4(1-G_u)(G_{uu}u_1^2+2G_{ux_1}u_1+G_{x_1x_1})\sum_{2\leq i\leq n}\frac{w_i^2}{w_1^2}u_{ii}\notag\\
=:&\sum_{2\leq i\leq n}K_iu_{ii},
\end{align}
where
\begin{align}\label{3Kia}
K_i=&-2(1-G_u)(h'\frac{Au_1}{w_1^3}-2\frac{D}{w_1^2})w_i^2
+2(1-G_u)g'A\frac{v^2+1}{w_1^2}w_i\gamma^i\notag\\
&-4(1-G_u)\frac{v^2+1}{w_1^2}(G_{ux_i}u_1w_1+\sum_{1\leq k\leq n} w_kG_{x_kx_i})w_i\notag\\
&+\big[-2h'\frac{Au_1^2}{v^2w_1^2}+4G_{uu}\frac{u_1^3}{v^2w_1}-2g'\gamma^1\frac{Au_1}{v^2w_1}+4G_{ux_1}\frac{u_1^2}{v^2w_1}
+4\sum_{1\leq k\leq n} w_kG_{ux_k}\frac{u_1^2}{v^2w_1^2}\notag\\&
+\sum_{1\leq k\leq n} w_kG_{x_kx_1}\frac{u_1}{v^2w_1^2}\big]w_i^2
+2g'\frac{Au_1}{w_1}w_i\gamma^i-4u_1^2w_iG_{ux_i}
-\frac{4u_1}{w_1}\sum_{1\leq k\leq n}w_kG_{x_kx_i}w_i\notag\\&
-4(1-G_u)v^2G_{x_ix_i}-4(1-G_u)(G_{uu}u_1^2+2G_{ux_1}u_1+G_{x_1x_1})\frac{w_i^2}{w_1^2}\notag\\
&+4(1-G_u)\frac{v^2+1}{w_1}(G_{ux_i}u_1+G_{x_1x_i})w_i.
\end{align}
It follows that
\begin{align}\label{3Kib}
|K_i|\leq&C_4A,\quad i=2,\ldots,n.
\end{align}

We write other terms as $J_2$, then
\begin{align}\label{3J2a}
J_2=&Q_3-h'fvu_1^2+\frac{2fG_{uu}}{A}vu_1^3w_1+\frac{h'^2}{2}\frac{Au_1^2}{w_1^2}\notag\\
+&\frac{h'^2}{1-G_u}\frac{Au_1^3}{v^2w_1}
+\frac{c^{11}}{2}g'^2\frac{u_1^2}{w_1^2}A
+\frac{c^{11}g'^2}{1-G_u}\frac{Au_1}{w_1}+O(u_1^2).
\end{align}

So by the choice of $\mu_0$ and the formula on $Q_3$ in  \eqref{3Q3b},\eqref{3A} and \eqref{3.1}, we get the following estimate on $J_2$,
\begin{align}\label{3J2b}
J_2\geq&-(1+c^{11}\alpha_0^2) u_1^2\log|Dw|^2+\frac{h'^2}{2}\frac{Au_1^2}{w_1^2}+\frac{h'^2}{1-G_u}\frac{Au_1^3}{v^2w_1}\notag\\&
+\frac{c^{11}}{2}g'^2\frac{u_1^2}{w_1^2}A
+\frac{c^{11}g'^2}{1-G_u}\frac{Au_1}{w_1}-C_5u_1^2\notag\\
\geq&\frac{1}{4}(1+c^{11}\alpha_0^2) u_1^2\log|Dw|^2-C_6u_1^2,
\end{align}
where  $ C_4, C_5, C_6 $ are  positive constants which only depend on $n, \Omega, \mu_0, M_0, L_1, L_2$.

{\bf Step 3:}
In this step, we concentrate on $J_1$. We first treat the terms $J_{11}$ and $J_{12}$  and obtain the formula \eqref{3aijvarphiijd},
then we complete the proof of  Theorem~\ref{Thm1.1} through  Lemma~\ref{Lem4.5}.

 By \eqref{3uii=2}, we have
\begin{align}\label{3u22}
\begin{split}
u_{22}
=&-\frac{1}{d_2}\sum_{3\leq i\leq n}d_iu_{ii}
+\frac{1}{d_2}[f v^3+\frac{h'}{2(1-G_u)}\frac{Au_1}{w_1}-\frac{D}{1-G_u}],
\end{split}
\end{align}

We first treat the term $J_{11}$:
using \eqref{3u22} to simplify \eqref{3J11a}, we get
\begin{align}\label{3J11b}
\begin{split}
J_{11}
=&2(1-G_u)^2\bigg\{d_2e_2u_{22}^2+\sum_{3\leq i\leq n}d_ie_iu_{ii}^2+2\frac{w_2^2}{w_1^2}u_{22}\sum_{3\leq j\leq n}\frac{w_j^2}{w_1^2}u_{jj}\\&
+2\sum_{3\leq i< j\leq n}\frac{w_i^2w_j^2}{w_1^4}u_{ii}u_{jj}
\bigg\}\\
=&2(1-G_u)^2\bigg\{\frac{e_2}{d_2}\big[-\sum_{3\leq i\leq n}d_iu_{ii}+f v^3+\frac{h'}{2(1-G_u)}\frac{Au_1}{w_1}-\frac{D}{1-G_u}\big]^2
+\sum_{3\leq i\leq n}d_ie_iu_{ii}^2\\&+\frac{2w_2^2}{d_2w_1^2}\big[-\sum_{3\leq i\leq n}d_iu_{ii}+f v^3
+\frac{h'}{2(1-G_u)}\frac{Au_1}{w_1}-\frac{D}{1-G_u}\big]\sum_{3\leq j\leq n}\frac{w_j^2}{w_1^2}u_{jj}\\&
+2\sum_{3\leq i< j\leq n}\frac{w_i^2w_j^2}{w_1^4}u_{ii}u_{jj}
\bigg\}\\
=&\frac{2(1-G_u)^2}{d_2}\bigg\{\sum_{3\leq i\leq n}\big[e_2d_i^2+e_id_id_2-2\frac{w_2^2w_i^2}{w_1^4}d_i\big]u_{ii}^2\\&
+2\sum_{3\leq i< j\leq n}\big[e_2d_id_j-\frac{w_2^2w_i^2}{w_1^4}d_j-\frac{w_2^2w_j^2}{w_1^4}d_i+\frac{w_i^2w_j^2}{w_1^4}\big]u_{ii}u_{jj}\\&
-2e_2\big[f v^3+\frac{h'}{2(1-G_u)}\frac{Au_1}{w_1}-\frac{D}{1-G_u}\big]\sum_{3\leq i\leq n}d_iu_{ii}\\&
+2\frac{w_2^2}{w_1^4}\big[f v^3+\frac{h'}{2(1-G_u)}\frac{Au_1}{w_1}-\frac{D}{1-G_u}\big]\sum_{3\leq i\leq n}w_i^2u_{ii}\\&
+e_2\big[f v^3+\frac{h'}{2(1-G_u)}\frac{Au_1}{w_1}-\frac{D}{1-G_u}\big]^2\bigg\}.
\end{split}
\end{align}

 We can rewrite it as the following
\begin{align}\label{3J11bb}
\begin{split}
J_{11}
  =&:\frac{2(1-G_u)^2}{d_2}\bigg\{\sum_{3\leq i\leq n}b_{ii}u_{ii}^2+2\sum_{3\leq i< j\leq n}b_{ij}u_{ii}u_{jj}\\&-2e_2\big[f v^3+\frac{h'}{2(1-G_u)}\frac{Au_1}{w_1}-\frac{D}{1-G_u}\big]\sum_{3\leq i\leq n}d_iu_{ii}\\&
+2\frac{w_2^2}{w_1^4}\big[f v^3+\frac{h'}{2(1-G_u)}\frac{Au_1}{w_1}-\frac{D}{1-G_u}\big]\sum_{3\leq i\leq n}w_i^2u_{ii}\\&
+e_2\big[f v^3+\frac{h'}{2(1-G_u)}\frac{Au_1}{w_1}-\frac{D}{1-G_u}\big]^2\bigg\},
\end{split}
\end{align}
where
\begin{align}\label{3bija}
\begin{split}
b_{ii}=&e_2d_i^2+e_id_id_2-2\frac{w_2^2w_i^2}{w_1^4}d_i\\
=&:2u_1^4+A_{1i}u_1^2+A_{2i},\quad i\geq 3\\
b_{ij}=&e_2d_id_j-\frac{w_2^2w_i^2}{w_1^4}d_j-\frac{w_2^2w_j^2}{w_1^4}d_i+\frac{w_i^2w_j^2}{w_1^4}\\
=&:u_1^4+G_{ij}u_1^2+\hat{G}_{ij},
\qquad\qquad\qquad\qquad i\neq j,\,i,j\geq 3,
\end{split}
\end{align}
and
\begin{align}\label{3A1i}
\begin{split}
A_{1i}=&4+(w_2^2+w_i^2)\frac{u_1^2}{w_1^2}, \\
A_{2i}=&2+(3w_2^2+5w_i^2)\frac{u_1^2}{w_1^2}
+w_i^2(w_2^2+w_i^2)\frac{u_1^2}{w_1^4}+\frac{2w_2^2+4w_i^2}{w_1^2}+\frac{2w_i^2(w_2^2+w_i^2)}{w_1^4}, \\
G_{ij}=&2+w_2^2\frac{u_1^2}{w_1^2}, \\
\hat{G}_{ij}=&1+(2w_2^2+w_i^2+w_j^2)\frac{u_1^2}{w_1^2}
+\frac{w_2^2+w_i^2+w_j^2}{w_1^2}+\frac{2w_i^2w_j^2}{w_1^4}-\frac{w_2^2w_i^2w_j^2}{w_1^6}.
\end{split}
\end{align}

Now we simplify the terms in $J_{12}$: by  \eqref{3u22},  we can rewrite \eqref{3J12a} as
\begin{align}\label{3J12b}
\begin{split}
J_{12}=&K_2u_{22}+\sum_{3\leq i\leq n}K_iu_{ii}\\
=&\frac{K_2}{d_2}\big[-\sum_{3\leq i\leq n}d_iu_{ii}
+f v^3+\frac{h'}{2(1-G_u)}\frac{Au_1}{w_1}-\frac{D}{1-G_u}\big]+\sum_{3\leq i\leq n}K_iu_{ii}\\
=&\sum_{3\leq i\leq n}[K_i-\frac{K_2d_i}{d_2}]u_{ii}
+\frac{K_2}{d_2}[f v^3+\frac{h'}{2(1-G_u)}\frac{Au_1}{w_1}-\frac{D}{1-G_u}].
\end{split}
\end{align}
Using  \eqref{3J11b} and \eqref{3J12b} to treat \eqref{3J1a},  we have
\begin{align}\label{3J1b}
\begin{split}
J_{1}=&\frac{2(1-G_u)^2}{Ad_2}\bigg\{\sum_{3\leq i\leq n}b_{ii}u_{ii}^2+2\sum_{3\leq i< j\leq n}b_{ij}u_{ii}u_{jj}\\&
+\sum_{3\leq i\leq n}\tilde{K}_iu_{ii}\big]+ R,
\end{split}
\end{align}
where
\begin{align}\label{3tildeKia}
\tilde{K}_i=&-2e_2\big[f v^3+\frac{h'}{2(1-G_u)}\frac{Au_1}{w_1}-\frac{D}{1-G_u}\big]d_i\notag\\&
+2\frac{w_2^2}{w_1^4}\big[f v^3+\frac{h'}{2(1-G_u)}\frac{Au_1}{w_1}-\frac{D}{1-G_u}\big]w_i^2+
K_id_2-K_2d_i.
\end{align}
We also have let
\begin{align*}
\begin{split}
R=&\frac{2e_2(1-G_u)^2}{Ad_2}\big[f v^3+\frac{h'}{2(1-G_u)}\frac{Au_1}{w_1}-\frac{D}{1-G_u}\big]^2\\
+&\frac{K_2}{Ad_2}[f v^3+\frac{h'}{2(1-G_u)}\frac{Au_1}{w_1}-\frac{D}{1-G_u}].
\end{split}
\end{align*}

For $\tilde{K}_i$ and $R$, using the formulas on $D, K_i$ in \eqref{3Db}, \eqref{3Kib}; the formula of $A$ in \eqref{3A}; $ e_i; d_i$  in \eqref{3ei}-\eqref{3di},
and  $h(u), g(d)$ in \eqref{3.1}, we have the following estimates
\begin{align}
|\tilde{K}_i|\leq& C_7u_1^5,\quad i=3,\ldots,n;\label{3tildeKib}\\
|R|\leq &C_8\frac{u_1^2}{\log u_1}.\label{3R}
\end{align}

Now we use  Lemma~\ref{Lem4.5}, if there is a sufficiently large positive constant $C_9$ such that
\begin{align}
|Du|(x_0) \ge C_9,\label{3C9}
\end{align}
 then we have
\begin{align}\label{3J1c}
J_{1}
\geq &\frac{2}{Ad_2}({-C_{10} u_1^6}),\notag\\
\ge & -C_{11} u_1^2           ,
\end{align}
where we use the formulas  $d_2$ in \eqref{3di} and $A$ in \eqref{3A}.

Using the estimates on $J_1$ in \eqref{3J1c} and  $J_2$ in \eqref{3J2b}, from  \eqref{3aijvarphiijc} we obtain
\begin{align}\label{3aijvarphiijd}
0\geq &\sum_{1\leq i,j\leq n}a^{ij}\varphi_{ij}\notag\\
 \geq&\frac{1}{4}(1+c^{11}\alpha_0^2) u_1^2\log|Dw|^2-C_{12}u_1^2\notag\\
 \geq&\frac{1}{4} u_1^2\log|Dw|^2-C_{12}u_1^2,
\end{align}

 There exists a positive constant
$C_{13}$ such that
 \begin{align}\label{3C13}
 |Du|(x_0)\leq C_{13}.
 \end{align}

 So from  Case 1, Case 2,  and \eqref{3C13}, we have
 $$|Du|(x_0)\leq C_{14}, \quad \quad x_0\in\Omega_{\mu_0}\bigcup\partial\Omega.$$

Since $\varphi(x)\leq\varphi(x_0),\quad \text{for} \quad x\in \Omega_{\mu_0}$, there exists $M_2$ such that
\begin{align}\label{3M2}
|Du|(x)\leq M_2, \quad in\quad\Omega_{\mu_0}\bigcup\partial\Omega,
\end{align}
where~$M_2$ depends only on ~$n, \Omega, \mu_0, M_0,  L_1, L_2$.

So at last we get the following estimate
$$\sup_{\overline\Omega_{\mu_0}}|Du|\leq \max\{M_1, M_2\},$$
where the positive constant  ~$ M_1$ depends only on $n, \mu_0, M_0, L_1$; and $ M_2$ depends only on ~$n, \Omega, \mu_0, M_0, L_1, L_2$.

 So  we  complete the proof of Theorem~\ref{Thm1.1}.\qed

   Now we
 prove the main Lemma~\ref{Lem4.5}  which was used to estimate  $J_1$ defined in  \eqref{3J1c}.

\begin{Lem}\label{Lem4.5}
We define $(b_{ij})$  as in \eqref{3bija}; $d_i ,e_i$ defined as in \eqref{3di}-\eqref{3ei}; $A_{1i}, A_{2i}, G_{ij}, \hat{G}_{ij}$ defined as in \eqref{3A1i}.
We study the following quadratic form
\begin{align}\label{4.6}
\begin{split}
Q(x_3,x_4,\ldots,x_n)=&\sum_{3\leq i\leq n}b_{ii}x_i^2
+2\sum_{3\leq i< j\leq n}b_{ij}x_i x_j
+\sum_{3\leq i\leq n}\tilde{K}_i x_i,
\end{split}
\end{align}
where $\tilde{K}_i$ defined in  \eqref{3tildeKia} and we have the estimate  \eqref{3tildeKib} for $\tilde{K}_i$.
Then there exists a sufficiently large positive constant $C_{15}$ which depends only on  $n, \Omega, \mu_0, M_0, L_1, L_2,$ such that if
\begin{align}\label{4C15}
|Du|(x_0)=u_1(x_0) \ge C_{15},
\end{align}
then  the followings hold.

(I): The matrix $(b_{ij})$ is positive definite since the matrix  ~$(b^1_{ij})= (1+\delta_{ij})$ is positive definite.\par

(II): We have
\begin{align}\label{4.7}
\begin{split}
Q(x_3,x_4,\ldots,x_n)
\geq &-C_{16} u_1^6,
\end{split}
\end{align}
where positive constant $C_{16}$ also depends only on  $n, \Omega, \mu_0, M_0, L_1, L_2.$

\end{Lem}
{\bf Proof:} Let $$ B=(b_{ij})=B_1+B_2, B_1=u_1^4(b^1_{ij}), B_2=(O(u^2_1)\delta_{ij}).$$
We first  prove  (I):
\begin{align}\label{4.8}
\begin{split}
\sigma_k(B)=&\sigma_k(B_1+B_2)\\
=&\sigma_k(B_1)+\sigma_k(B_1, B_1, \ldots, B_1, B_2)\\&+\cdots+\sigma_k(B_1, B_2,\ldots, B_2, B_2)+\sigma_k(B_2)\\
=&u_1^{4k}\sigma_k(b^1_{ij})+O(u_1^{4k-2}),
\end{split}
\end{align}
so if $u_1$ is sufficiently large, then $\sigma_k(B)>0\Longleftrightarrow\sigma_k(b^1_{ij})>0$.

Now we prove   (II):
Since  ~$B_1=u_1^4(b^1_{ij})_{3\leq i,j\leq n}$  is positive definite, from the argument in (I), we get
\begin{align}\label{4.9}
B^{-1}=(B_1+B_2)^{-1}=B_1^{-1}(I+B_1^{-1}B_2)^{-1}=u_1^{-4}(b^1_{ij})^{-1}\big(1+o(1)\big).
\end{align}
Then we have
\begin{align}\label{4.10}
\begin{split}
(b^1_{ij})^{-1}=&\left(\begin{array}{cccc}2&1&\cdots &1\\
1&2&\cdots &1\\
\vdots & \vdots & \vdots &\vdots \\
1&1&\cdots& 2
\end{array}\right)^{-1}
=\frac{1}{n-1}\left(\begin{array}{cccc}
n-2&-1&\cdots & -1\\
-1&n-2&\cdots & -1\\
\vdots & \vdots & \vdots &\vdots \\
-1&-1&\cdots & n-2
\end{array}\right).
\end{split}
\end{align}

Now we solve the following linear  algebra equation
\begin{align}\label{4.11}
\frac{\partial Q}{\partial x_k}=0,\quad k=3,4,\ldots,n.
\end{align}
 We assume $(\bar{x}_{3},\bar{x}_{4},\ldots,\bar{x}_{n})$ is the extreme point of the quadratic form $Q(x_3,x_4,\ldots,x_n)$.
From the definition of $b_{ij}, \tilde{K}_i$ in ~\eqref{3bija}, ~\eqref{3tildeKia} and the estimate for $\tilde{K}_i$ in ~\eqref{3tildeKib}, using the formulas \eqref{4.9} and \eqref{4.10}, it follows that
\begin{align}\label{4.12}
\begin{split}
\left(\begin{array}{c}
\bar{x}_{3}\\
\bar{x}_{4}\\
\vdots  \\
\bar{x}_{n}
\end{array}\right)=& O(u_1^5)B^{-1}
\left(\begin{array}{c}
1\\
1\\
\vdots\\
1
\end{array}\right)
=O(u_1)
\left(\begin{array}{c}
1\\
1\\
\vdots\\
1
\end{array}\right).
\end{split}
\end{align}

It follows that we have the following minimum of the quadratic $Q$,
\begin{align}\label{4.14}
\begin{split}
Q(\bar{x}_3, \bar{x}_4,\ldots, \bar{x}_n)
\geq&
-C_{17}u_1^6.
\end{split}
\end{align}

 In this computation, the bounds in the coefficient on $ O(u_1^5), O(u_1)$ depend only on $n, \Omega,  M_0, \mu_0,  L_1, L_2$. Thus we complete this proof.\qed

\section{ The proof of Theorem \ref{Thm1.2}  }

In this section we first prove the Theorem~\ref{Thm1.2}.

  In the proof of the existence theorem for the Neumann boundary value problem, we need the a priori estimates. For the $C^0$ estimates we use the methods introduced by ~Concus-Finn\cite{CF74} and ~Spruck\cite{Sp75}. As in  Simon-Spruck\cite{SS76}, we use the continuity method to complete the proof of theorem ~\ref{Thm1.2}

{\em Proof of Theorem~\ref{Thm1.2}:}

We consider the following family on the mean curvature equation with Neumann boundary value problem:
\begin{align}
\texttt{div}(\frac{Du}{\sqrt{1+|Du|^2}}) =&u   \quad\text{in}\quad\Omega,\label{5.1} \\
              \frac{\partial u}{\partial \gamma} =&\tau\psi(x)  \quad\text{on}\quad \partial \Omega,\label{5.2}
\end{align}
where $\tau\in[0,1]$.

For
$\tau=0$, then  $u=0$ is the unique solution. And we need to find the solution for
$\tau=1.$
By the standard existence  theorem \cite{Ur73, LU68}, as in Simon-Spruck \cite{SS76}, if we can get the a priori estimates for the $C^2(\overline\Omega)$ solution of the equation  \eqref{5.1} and \eqref{5.2}
\begin{align}
\sup_{\Omega}|u|\leq& \hat{K}_1,\label{5.3}\\
\sup_{\Omega}|Du|\leq& \hat{K}_2,\label{5.4}
\end{align}
where ~$\hat{K}_1, \hat{K}_2$  are independent of~$\tau$.
Then we can get the existence theorem. From the interior gradient estimates and our boundary gradient estimates, we only need get the $C^0$ estimates for the solution $u$ in  \eqref{5.1} and \eqref{5.2}.

In the paper by  Spruck\cite{Sp75}, he used the comparison theorem developed by Concus-Finn\cite{CF74} to get the $C^0$ estimates for the mean curvature equation with prescribed contact angle boundary value problem. In our case, his proof is still true, so we complete the proof of Theorem \ref{Thm1.2}.\qed

\begin{Rem}\label{4Rem2}
In X.N. Ma  and J.J. Xu \cite{MX14}, we have generalized the boundary gradient estimates to Hessian equation and the higher order curvature equation with Neumann boundary value and the capillary boundary value problem. And  we can reduce the condition on $\psi \in C^3$ to
$\psi \in  C^2$, then we could get the similar results  on Theorem 1.1 and Theorem 1.3.
\end{Rem}

      {\bf
Acknowledgments:}
The authors would like to
thank Qiu Guohuan for
 pointing
a flaw in our first version of this paper.

\end{document}